\documentclass[12pt]{article}
\usepackage{amssymb,amsmath,amsthm}

\def\N{{\mathbb N}}

\def\sobre#1#2{\lower 1ex \hbox{ $#1 \atop #2 $ } }

\newtheorem{theorem}{Theorem}[section]
\newtheorem{lemma}[theorem]{Lemma}

\renewenvironment{proof}[1]
{\noindent{\bf Proof.}\hspace{0.1cm} #1} {$\;\blacksquare$\newline}

\makeatletter
\@addtoreset{equation}{section}

\makeatother

\newcommand{\XC}{(X^{[k]}_n)_{n \in {\mathbb Z}}}

\newcommand{\X}{(X_n)_{n \in {\mathbb Z}}}

\newcommand{\XKn}{(X^{(k)}_n)_{n \in {\mathbb Z}}}
\newcommand{\pro}{{\mathbb P}}
\newcommand{\esp}{{\mathbb E}}
\newcommand{\var}{{\mathbb Var}}
\newcommand{\covar}{{\mathbb Cov}}
\newcommand{\one}{{\mathbf 1}}
\newcommand{\infd}{\underline{\delta}}
\newcommand{\ZK}{Z}
\newcommand{\ZM}{\tilde{Z}}
\newcommand{\Oun}{C}

\font\gfontb=cmbx10 scaled\magstep2
\newcommand{\gun}{\hbox{\gfontb 1}}

\title{Bootstrap Central Limit Theorem for Chains of Infinite Order
via Markov Approximations
\footnotetext{This work was partially supported by USP-Cofecub
agreement(project {\sl Processus stochastiques et m\'ecanique
statistique math\'ematique}), Brazil-France mathematics agreement,
CNPq (grant 301301-79) and FAPESP's Project {\sl Stochastic behavior,
critical phenomena and rhythmic pattern identification in natural
languages} (project number 03/09930-9).\\
{\em AMS 2000 subject classification. Primary}\/ 62F40,
60F05; {\em secondary} 60G10, 62G09.\\
{\em Keywords and phrases}:
bootstrap, chains of infinite order, canonical Markov approximations,
central limit theorem.\\
{\em Running head:} Bootstrap via Markov approximations.}
}

\author{P. Collet\\
Centre de Physique Th\'eorique\\
CNRS UMR 7644\\
Ecole Polytechnique\\
F-91128 Palaiseau Cedex, France \\
e-mail: \tt{collet@cpht.polytechnique.fr}\\
\and D. Duarte\\
Instituto de Matem\'atica e Estat\'\i stica\\
Universidade Federal de Goi\'as\\
Campus II, Goi\^ania, Brasil \\
e-mail: \tt{denise@mat.ufg.br}\\
\and A. Galves\\
Instituto de Matem\'atica e Estat\'\i stica\\
Universidade de S\~ao Paulo\\
Rua do Mat\~ao, 1010 \\
05508-090 S\~ao Paulo SP, Brasil \\
e-mail: \tt{galves@ime.usp.br}\\
}
\date{}

\begin{document}

\maketitle

\section*{Abstract}

We present a new approach to the bootstrap for chains of infinite
order taking values on a finite alphabet. It is based on a sequential
Bootstrap Central Limit Theorem for the sequence of canonical Markov
approximations of the chain of infinite order.  Combined with previous
results on the rate of approximation this leads to a Central Limit
Theorem for the bootstrapped estimator of the sample mean which is the
main result of this paper.

\section{Introduction.}\label{intro} 

In this paper we introduce a new procedure of bootstrap resampling for
chains on a finite alphabet whose transition probabilities depend on
the whole past. This resampling uses the excursions of the chain
between successive occurrences of the initial string of $k$ symbols as
building blocks for the bootstrap sample. The bootstrap sample is
obtained by concatenating randomly chosen blocks. These blocks are
chosen uniformly and independently among the first $m_k$ excursion
blocks. For chains which lose memory exponentially fast we prove a
Central Limit Theorem for the empirical mean of the bootstrap sample,
when the length $k$ of the initial reference string as well as the
number of excursion blocks $m_k$ diverge with a suitable relation
between them. This is the main result of the article.

The idea behind our procedure is that a typical large sample of the
chain of infinite order behaves essentially as a sample of a Markov
chain of order $k$ suitably chosen. The Markov property of the
approximating chain implies that the successive excursion blocks are
independent and identically distributed. This makes it possible to
construct the bootstrap sample by simply concatenating randomly chosen
blocks, exactly as proposed in the original paper by Efron
(1979)\nocite{efron} for the case of i.i.d. random variables.

This idea has already been exploited in the case of Markov chains in Athreya
and Fuh (1992)\nocite{af}. For chains of infinite order with different types
of mixing conditions, different approaches to the bootstrap have been proposed
in the papers by Calrstein (1986)\nocite{carl} and K\"unsch
(1989)\nocite{kunsch} and thoroughly studied in the recent literature, see for
example Liu and Singh (1992)\nocite{ls}, Shao and Yu (1993)\nocite{sy},
Naik-Ninbalkar and Rajarsh (1994)\nocite{N}, B\"uhlmann (1994)\nocite{bul} and
Peligrad (1998)\nocite{mp}.

Chains of infinite order seem to have been first studied by Onicescu and Mihoc
(1935a) \nocite{oni} who called them \emph{cha\^{\i}nes \`a liaisons
  compl\`etes}.  Their study was soon taken up by Doeblin and Fortet (1937)
\nocite{doefor37} who proved the first results on speed of convergence towards
the invariant measure.  The name chains of infinite order was coined by Harris
(1955) \nocite{har55}.  Our proof is based on the upper bound on the rate of
approximation of the chain of infinite order by the sequence of canonical
Markov approximations presented in Fern\'andez and Galves (2002)\nocite{fg}.
We also use the $\varphi$-mixing property of the chain of infinite order
proven in Bressaud, Fern\'andez and Galves (1999)\nocite{bfg2}.  We refer the
reader to Iosifescu and Grigorescu (1990) \nocite{iosgri90} for a complete
survey, and to Fern\'andez, Ferrari and Galves, 2001\nocite{ffg}) for an
elementary presentation of the subject from a constructive point of view.

The rest of the paper is organized as follows. In section \ref{nota}
we introduce the notation and the definitions and state the main
results. In section \ref{preliminary} we collect together a few
technical results which will be used in the proof of the theorems. In
section \ref{proofbootm} we prove a central limit theorem for the
sequence of canonical approximating Markov chains. Finally in section
\ref{proofboot} we prove the main result which is a bootstrap central
limit theorem for the empirical mean of a chain of infinite order.

\section{Notations, definitions and statement of the main result.}\label{nota}
Let $\X$ be a stationary process taking values on a finite alphabet
$A$.  We will use the shorthand notation
$$
p(x_0|x_{-1},x_{-2},\ldots)=\pro\big(X_0=x_0\big|
X_{-1}=x_{-1},X_{-2}=x_{-2},\ldots\big)
$$  
to denote the regular version of the conditional probability of the
process. To avoid long formulas, whenever convenient, we will use the
notation $a_{0,l}$ to denote the sequence $(a_0,\ldots, a_l)$ of
elements of $A$. We also use the notation $\{X_{n, n+l}= a_{0,l}\}$ to
denote the cylinder set
$$\{X_n=a_0, \ldots, X_{n+l}=a_l\}\;.$$
Following Harris (1955), we call this process a chain of infinite order.

We assume that $\X$ satisfies the following hypotheses.

{$\bf H_1$}
\begin{equation}\label{p3}
 \min_{a \in A}\inf_{ (\ldots, x_{-2}, x_{-1}) \in  {\cal{A}}_a } p(a
|x_{-1}, x_{-2}, \ldots) =\delta >0\;,
\end{equation}
where
\noindent ${\cal{A}}_a=
\{ (\ldots, x_{-2}, x_{-1}):p(a| x_{-1}, x_{-2}, \ldots) >0\}$.

{$\bf H_2$}
$$
c=-\limsup_{l \to \infty}\frac{1}{l}\log \beta_l>0\; ,
$$
where
\[
  \beta_l =
 \sup_{ \stackrel{x_i=y_i} {i=-l,\ldots,0}}
\big|p(x_0 |x_{-1},x_{-2},\ldots) - p(y_0 |
y_{-1},y_{-2},\ldots)\big|\; .
\]

Let $f:A^r \rightarrow {\mathbb R}$ be a real observable of the chain,
where $r$ is a fixed positive integer and denote\ \ \ \ \ \ 
         
$$
\mu=\esp\big(f(X_1,\ldots,X_r)\big)\;,
$$ the average value of the observable $f$.  We are interested in the
fluctuations of an estimator of $\mu$. To simplify the presentation
we can assume without loss of generality that $r=1$, namely the
cylinder function $f$ through which we observe the chain depends only
on one coordinate.

To avoid uninteresting pathologies we will assume that the following
third hypothesis holds

{$\bf H_3$}
\[
\sigma^2=\var\left(f(X_0)\right)+
2\sum_{j=1}^{+\infty}\covar\left(f(X_0), f(X_j)\right) >0\; .
\]

We recall that hypotheses {$\bf H_1$} and {$\bf H_2$} imply that the
chain $\X$ is exponentially $\varphi$-mixing (cf. Bressaud,
Fern\'andez and Galves (1999)\nocite{bfg2}). This last property imply
that the series defining $\sigma^2$ is convergent (cf. for instance
Theorem 19.1 in Billingsley 1999). However it is well known that this
does not imply that $\sigma$ is strictly positive. 

Our bootstrap procedure is defined as follows. For any positive
integer $k$, the sequence $\big(R_j(k)\big)_{j\in\N}$ of return times
of the first string of length $k$ is defined by
\[
R_{i+1}(k)=\inf\bigg\{n>R_{i}(k)\;:\;\big(X_n,\ldots,X_{n+k-1}\big)=
\big(X_1,\ldots,X_k\big)\bigg\}\;,
\]
with $R_0(k)=1$.

Let $\xi_i(k)$ be the  block of values of the chain 
from $R_{i-1}(k)$ up to $R_{i}(k)-1$,  namely
\begin{equation}\label{block}
\xi_{i}(k) \,=\,\big(X_{R_{i-1}(k)},\ldots, X_{R_{i}(k)-1}\big)\, .
\end{equation}
We will make a uniform i.i.d. selection of the first $m$ blocks
$\xi_{1}(k),\ldots,\xi_{m}(k)$ to construct a bootstrap sample of the
chain.  We will take $m=m_k$ as a diverging function of $k$ to be
fixed latter. This leads naturally to the construction of a sequence
of bootstrap samples indexed by $k$.

The formal definition is the following. For every $k$, let
$I_1(k),\ldots,I_{m_k}(k)$ be $m_k$ independent random variables with
uniform distribution in the set $\{1,\ldots,m_k\}$. The bootstrap
blocks are defined as
\[
\xi^{\ast}_{l}(k)=\xi_{I_l(k)}(k)\, ,
\]
for $l=1,\ldots,m_k$. The bootstrap sample $X^*_1(k),\ldots
X^*_{R^*_{m_k}(k)}(k)$ is constructed by concatenating the bootstrap
blocks $ \xi^{\ast}_{1}(k),\dots,\xi^{\ast}_{m(k)}(k)$.  We observe
that the return times of the bootstrap sample assume the values
$ R^*_0(k)=1 $
and for $l=1,\ldots,m_k$
$$
R^{*}_{l}(k)=R^{*}_{l-1}(k)+R_{I_l(k)+1}(k)-R_{I_l(k)}(k)\;.
$$

We consider the following sequence of estimators for $\mu$
\begin{equation}\label{xbarb}
\hat{\mu}_{k} = \frac{1}{R_{m_k}(k)-1}\sum^{R_{m_k}(k)-1}_{n=1} f(X_n)\; .
\end{equation}
Its bootstrap counterpart is given by
\begin{equation} 
\label{xbarboot}
\mu^*_{k}=\frac{1}{R^*_{m_k}(k)-1}\sum^{R^*_{m_k}(k)-1}_{n=1} f(X^*_n(k))\;.
\end{equation}

Let
$$
\sigma_k^*= \sqrt{\frac{\var\left(\sum_{n=1}^{R_{m_k}^*(k)-1}
\big(f(X^*_n(k))-\hat{\mu}_k\big)\;
      \bigg|\;X_1,\ldots,X_{R_{m_k}-1}\right)}{R_{m_k}^*(k)-1}} \;,
$$
where $\var$ denotes the variance.  Observe that $\sigma_k^*$ is a
function of the sample $X_1,\ldots,X_{R_{m_k}(k)}$ and therefore the
above variance is taken with respect to the independent random
variables $I_1(k),\ldots, I_{m_k}(k)$.

In the statement of our theorems the number of blocks used in the
bootstrap sample is 
\[
m_k(\alpha)=[e^{\alpha k}]\; ,
\]
where $\alpha$ is a positive real number to be suitably chosen latter
and $[\cdot]$ denotes the integer part. To simplify the notation we
will often write $m_k$ instead of $m_k(\alpha)$ 

\begin{theorem}\label{boot}
Let $\X$ be a chain of infinite order satisfying Hypotheses ${\bf
H_1}$, ${\bf H_2}$ and ${\bf H_3}$ and such that 
$c > 18\ln\left(1/\delta\right)$
, where $\delta$ and $c$ are the constants appearing in
${\bf H_1}$ and ${\bf H_2}$, respectively. Then, for any 
$\alpha \in \left(5 \ln\left(1/ \delta \right), c-\ln\left(1/
\delta \right) \right)$,
for $ m_{k}=[ e^{\alpha k}]$, and
for almost all realizations of the chain $\X$, we have
\begin{equation}
\frac{\sqrt{R^*_{m_k}(k)-1}}{\sigma_k^*}\big( \mu^*_{k}-
\hat{\mu}_{k} 
\big)\;
\stackrel{{\mathcal D}}{\longrightarrow} {\cal N}(0,1)\;,
\end{equation}
as $k$ tends to $+\infty$, where 
$\stackrel{{\mathcal D}}{\longrightarrow}$ 
denotes convergence in distribution and ${\cal N}(0,1)$ denotes the
standard normal distribution.
\end{theorem}

The proof of Theorem \ref{boot} is based on the following sequential
bootstrap procedure which is interesting by itself. Let $\XKn$,
$k=1,2,\ldots$ be a sequence of stationary irreducible aperiodic
Markov chains of order $k=1,2,\dots$, respectively, taking values in
the same finite alphabet $A$ with transition probabilities denoted by

$$
p^{(k)}(a|\; b_{-k, -1})=\pro(X_0=a|\; X_{-k, -1}^{(k)}=b_{-k,-1})\, .
$$
We may assume, without loss of generality, that the Markov chains
$\XKn$, for $k=1,2,\ldots$ are all defined on the same probability
space (cf. for instance \cite{ffg}).

We define
$$
\delta^{(k)}= 
\min_{a \in A}\inf_{ ( x_{-k}, \ldots , x_{-1})\in {\cal{A}}^{(k)}_a } 
p^{(k)}(a |x_{-1}, \ldots, x_{-k})\; 
$$
and
\begin{equation}\label{deltabarre}
\infd = \inf\{\delta^{(k)}:\; k\ge 1\}\; ,
\end{equation}
where ${\cal{A}}^{(k)}_a= \left\{ (x_{-k},\ldots, x_{-1}):p^{(k)}(a|
x_{-1}, \ldots, x_{-k}) >0\right\}$.

For each $k$ we define recursively the sequence of return times
$\big(R_j^{(k)}\big)_{j\in\N}$ by $R_0^{(k)}=1$, and for $i\ge 1$
\begin{equation}\label{recursi}
R_{i}^{(k)}=
\inf\bigg\{n>R_{i-1}^{(k)}\;:\;\big(X_n^{(k)},\ldots,X_{n+k-1}^{(k)}\big)=
\big(X_1^{(k)},\ldots,X_k^{(k)}\big)\bigg\}\;.
\end{equation}

Let $\xi^{(k)}_i$ be the block of values of the chain $\XKn$ from
$R_{i-1}^{(k)}$ up to $R_{i}^{(k)}-1$, namely
$$
\xi^{(k)}_i=\big(X^{(k)}_{R^{(k)}_{i-1}}, \ldots,
X^{(k)}_{R^{(k)}_{i}-1}\big)\;.
$$

We construct a bootstrap sample of the Markov chain $\XKn$ by
performing an i.i.d. selection of the blocks $\xi^{(k)}_l$. The formal
definition is the following. For every $k$, let
$I_1(k),\ldots,I_{m_k}(k)$ be $m_k$ independent random variables with
uniform distribution in the set $\{1,\ldots,m_k\}$.  The bootstrap
blocks are defined by
\[
\xi^{(k)*}_{l}=\xi^{(k)}_{I_l(k)}\, ,
\]
for $l=1,\ldots,m_k$. The bootstrap sample $X^{(k)*}_l,
l=1,\ldots,R^{(k)*}_{m_k}$ is constructed by concatenating the blocks
$\xi^{(k)*}_{1},\ldots, \xi^{(k)*}_{m_k}\;.$ We observe that the
return times of the bootstrap sample assume the values $R^{(k)*}_0=1$
and for $l=1,\ldots,m_k$
\[
R^{(k)*}_{l}=R^{(k)*}_{l-1}+R^{(k)}_{I_l(k)+1}-R^{(k)}_{I_l(k)}\;.
\]
We consider the following estimator for $\mu^{(k)}= \esp\big(
f(X^{(k)}_1)\big)$ 
\begin{equation}\label{xbarxk} \hat{\mu}^{(k)} =
  \frac{1}{R^{(k)}_{m_k}-1}\sum^{R^{(k)}_{m_k}-1}_{n=1}f( X^{(k)}_n)\;.
\end{equation}
Its bootstrap counterpart is given by
\begin{equation} 
\label{xbarbootk}
\mu^{(k)*}=\frac{1}{R^{(k)*}_{m_k}-1}\sum^{R^{(k)*}_{m_k}-1}_{n=1}f( X^{(k)*}_n)
\end{equation}

We define
\begin{equation}\label{sigmak}
\sigma^{(k)*}=
\sqrt{
\frac{\var \left(\sum_{n=1}^{R^{(k)*}_{m_k}-1}(f(X^{(k)*}_n)-\hat{\mu}^{(k)})\;
\bigg|\;X^{(k)}_1,\ldots,X^{(k)}_{R^{(k)}_{m_k}-1}\right)}
{R^{(k)*}_{m_k}-1}}
\end{equation}
Recall that, as before, this variance is with respect to the
independent random variables $I_1(k), \ldots, I_{m_k}(k)$.

\begin{theorem}
\label{bootmarkov} Let $\XKn$,
$k=1,2,\ldots$ be a sequence of stationary, irreducible, and aperiodic
Markov chains of order $k=1,2,\dots$, respectively, taking values in
the same finite alphabet $A$ and satisfying the following hypotheses 
\begin{equation}\label{deltauni}
\infd > 0\; ,
\end{equation}
where $\infd$ is defined in (\ref{deltabarre}), and
\begin{equation}\label{Ak}
\liminf_{k\to+\infty}\esp\left(\left(\sum_{n=1}^{R^{(k)}_{1}-1}
\left(f\left(X^{(k)}_n\right)-\mu^{(k)}\right)\right)^2\right)>0\; .
\end{equation}
If $\alpha > 5\ln\left(1/\infd\right)$ and $m_k= [e^{\alpha k}]$, then
for almost all realizations of the chains $\XKn$,  $k=1,2\ldots$, we
have
\[
\frac{\sqrt{R^{(k)*}_{m_k}-1}}
{\sigma^{(k)*}}\left(\mu^{(k)*}
-
\hat{\mu}^{(k)} \right)
\stackrel{{\mathcal D}}{\longrightarrow}
{\cal N}(0,1)\; ,
\]
as $k$ tends to $+\infty$.
\end{theorem}

\section{Preliminary results}\label{preliminary}

We first introduce some shorthand notation. We define
\[
{\ZK}^{(k)}_i=
\sum_{n=R^{(k)}_{i-1}}^{R^{(k)}_{i}-1}
\left(f\left(X^{(k)}_n\right)-\hat{\mu}^{(k)}\right)
\;,
\]
and its bootstrap version  is given by
\[{\ZK}^{(k)*}_i=
\sum_{n=R^{(k)*}_{i-1}}^{R^{(k)*}_{i}-1}
\left(f\left(X^{(k)*}_n\right)-\hat{\mu}^{(k)}\right)\;.
\]
Note that ${\ZK}^{(k)*}_i={\ZK}^{(k)}_{I_i(k)} $.

We use the shorthand $\esp^*(\,\cdot\,)$ to denote
$\esp\left(\,\cdot\,\big|X^{(k)}_1, \ldots,
X^{(k)}_{R_{m_k}^{(k)}}\right)$ and $\var^*(\,\cdot\,)$ to denote
$\var\left(\,\cdot \,\big|X^{(k)}_1, \ldots,
X^{(k)}_{R_{m_k}^{(k)}}\right)$. We recall that, in both cases, the
expectation is taken with respect to the sequence $I_i{(k)}$,
$i=1,\ldots,m_{k}$  of i.i.d. random variables uniformly distributed in
the set $\left\{1,\ldots,m_{k}\right\}$.

\begin{lemma}\label{expresboot}
The following equalities hold
\[
\esp^{*}\left({\ZK}_{1}^{(k)*}\right)=0 \;,
\]
and
\[
\var^{*}\left(\sum_{l=1}^{m_{k}}{\ZK}_{l}^{(k)*}\right)= \sum^{m_k}_{l=1}
\left({\ZK}^{(k)}_l\right)^2\;.
\]
\end{lemma}

\begin{proof}
By definition
\begin{equation}\label{espkstar}
\esp^{*}\left({\ZK}_{1}^{(k)*}\right)=
\sum_{n=1}^{m_{k}}{\ZK}_{n}^{(k)}\pro\left(I_1^{(k)}=n\right)=
\frac{1}{m_k}\sum_{n=1}^{m_{k}}{\ZK}_{n}^{(k)}=0\;.
\end{equation}
The second equality follows by a similar computation.
\end{proof}

It is convenient to introduce a new family of random variables
$\ZM^{(k)}_{i}$, where $i=1,\ldots ,m_k$, defined as follows
\begin{equation}\label{zm}
{\ZM}^{(k)}_i=
\sum_{n=R^{(k)}_{i-1}}^{R^{(k)}_{i}-1}
\left(f\left(X^{(k)}_n\right)-\mu^{(k)}\right)
\;.
\end{equation}
These random variables are not only identically distributed (as it was
already the case for $\big(\ZK^{(k)}_{l}\big)$), but also they are
independent and have zero mean. Moreover the following relation holds
\begin{equation}\label{dif}
\ZK^{(k)}_{l}=\ZM^{(k)}_{l}+\big(\mu^{(k)}-\hat{\mu}^{(k)}
\big)\big(R^{(k)}_{l}-R^{(k)}_{l-1}\big)
\end{equation}

We define $D^{(k)}_{l}=R^{(k)}_{l} -
R^{(k)}_{l-1}$ (recall that $ R^{(k)}_{0}=1$). Similarly, we define
$D^{(k)*}_{l}= R^{(k)*}_{l} - R^{(k)*}_{l-1}$.

\begin{lemma}\label{basic}
There is a positive constant $C$ independent of $k$ such that
\[
\big|{\ZK}_1^{(k)}\big| \le  C D_1^{(k)}\;,
\quad \mathrm{and}\quad 
\big|{\ZM}_1^{(k)}\big| \le  C D_1^{(k)}\;.
\]
\end{lemma}
 
\begin{proof}
This result follows immediately from the fact that the observable $f$
has finite range.
\end{proof} 

\begin{lemma}\label{suchthat}
There is a constant $C>0$ such that, for any $k\ge
1$, the following inequality holds
\[
\esp\left(\big(\hat{\mu}^{(k)}-\mu^{(k)}\big)^2\right)
\le C\, \frac{\esp\left(\big(D_1^{(k)}\big)^4\right)}{m_k}
\, .
\]
\end{lemma}
\begin{proof}
By definition we have
\[
\hat{\mu}^{(k)}-\mu^{(k)}=\frac{\sum_{l=1}^{m_k}{\ZM}_{l}^{(k)}}
{\sum_{l=1}^{m_k}{D}_{l}^{(k)}}\, 
\]
and therefore, using the Markov property and the stationarity of the
chain, we have
\begin{equation}\label{handside}
\esp\left(\big(\hat{\mu}^{(k)}-\mu^{(k)}\big)^2\right)=
\end{equation}
\begin{equation}
m_k\esp\left(\frac{\big({\ZM}_{1}^{(k)}\big)^2}
{\big(\sum_{l=1}^{m_k}{D}_{l}^{(k)}\big)^2}\right)
+
m_k(m_k-1)\esp\left(\frac{{\ZM}_{1}^{(k)}{\ZM}_{2}^{(k)}}
{\big(\sum_{l=1}^{m_k}{D}_{l}^{(k)}\big)^2}\right)\, .
\end{equation}
Since $\sum_{l=1}^{m_k}{D}_{l}^{(k)}>m_{k}$, and using Lemma
\ref{basic}, we conclude that the first term in the right hand side 
of expression (\ref{handside}) is bounded above by
\begin{equation}\label{premier}
C\;\frac{\esp\left(\big({D}_{1}^{(k)}\big)^2\right)}{m_{k}}
\end{equation}
where $C>0$ is a constant independent of $k$.

To obtain an upper bound for the second term on the right hand side of
expression (\ref{handside}), we first observe that for ${m_k}\ge 4$ we
have
\begin{equation}\label{compile}
\esp\left(\frac{{\ZM}_{1}^{(k)}{\ZM}_{2}^{(k)}}
{\big(\sum_{l=1}^{m_k}{D}_{l}^{(k)}\big)^2}\right)=
\esp\left(\frac{{\ZM}_{1}^{(k)}{\ZM}_{2}^{(k)}}
{\big(\sum_{l=3}^{m_k}{D}_{l}^{(k)}\big)^2}\right)
\end{equation}
\[
-\esp\left(\frac{{\ZM}_{1}^{(k)}{\ZM}_{2}^{(k)}
\big({D}_{1}^{(k)}+{D}_{2}^{(k)}\big)^2 }
{\big(\sum_{l=3}^{m_k}{D}_{l}^{(k)}\big)^2
\big(\sum_{l=1}^{m_k}{D}_{l}^{(k)}\big)^2}\right)
-2\esp\left(\frac{{\ZM}_{1}^{(k)}{\ZM}_{2}^{(k)}
\big({D}_{1}^{(k)}+{D}_{2}^{(k)}\big)}
{\big(\sum_{l=3}^{m_k}{D}_{l}^{(k)}\big)^2
\big(\sum_{l=1}^{m_k}{D}_{l}^{(k)}\big)^2}\right)\;.
\]

The independence of ${\ZM}_{1}^{(k)}$,  ${\ZM}_{2}^{(k)}$ and 
 $\sum_{l=3}^{m_k}{D}_{l}^{(k)}$ imply that

$$
\esp\left(\frac{{\ZM}_{1}^{(k)}{\ZM}_{2}^{(k)}}
{\big(\sum_{l=3}^{m_k}{D}_{l}^{(k)}\big)^2}\right)=0\;.
$$ 

Using again Lemma \ref{basic}, H\"older's inequality and
${D}_{l}^{(k)}\ge 1$, we deduce that the sum of the absolute values of
the two remaining terms of the right hand side of expression
(\ref{compile}) is bounded above by
\begin{equation}\label{deuxieme}
 C\, \frac{\esp\left(\big(D_1^{(k)}\big)^3\right)}{m_k^3}+
\frac{\esp\left(\big(D_1^{(k)}\big)^4\right)}{m_k^3}
 \,,
\end{equation}
where $C$ is a positive constant independent of $k$.
Since ${D}_{1}^{(k)}\ge 1 $, 
inequalities (\ref{premier}) and (\ref{deuxieme})  conclude the proof.
\end{proof}

\begin{lemma}\label{upprod1}
For any integer $k$ and any positive real
number $t$ the following inequality holds
\[
\pro\left(D_1^{(k)} >t\right) \le
\left(1-\infd^k\right)^{[t/k]}
\;.
\]
\end{lemma}

\begin{proof} 
We observe that
\[
\pro\left(D_1^{(k)} >t\right) \le
\pro\left(\bigcap_{j=1}^{[t/k] }
\left\{X_{jk+1, (j+1)k}^{(k)}\neq X_{1,k}^{(k)}\right\} \right)\; .
\]
Now we rewrite the right-hand side of the above inequality, by
conditioning on the values of the initial $k$ symbols
\[
\sum_{a_{1,k}}\pro\left(X_{1, k}^{(k)}=
a_{1,k}\right)\pro\left(\bigcap_{j=1}^{[t/k] }
\left\{X_{jk+1, (j+1)k}^{(k)}\neq a_{1,k}\right\} \, 
\bigg|\, X_{1, k}^{(k)}=a_{1,k}\right)\; .
\]
The second factor in the above sum can be rewritten as
\[
\left[1-\pro\left(X_{[t/k]k+1, ([t/k]+1)k}^{(k)}= a_{1,k}\,\bigg|
\bigcap_{j=1}^{[t/k]-1 }
\left\{X_{jk+1, (j+1)k}^{(k)}\neq a_{1,k}\right\}\bigcap
\left\{X_{1, k}^{(k)}=a_{1,k}\right\}\right) \right]
\]
\[
\times \;\pro\left(\bigcap_{j=1}^{[t/k]-1 }
\left\{X_{jk+1, (j+1)k}^{(k)}\neq a_{1,k}\right\}
\, \bigg| \, X_{1, k}^{(k)}=a_{1,k} \right)\, .
\]
Using  \ref{deltauni} this last expression can be
bounded above by
\[
\left(1-\infd^k\right)\pro\left(\bigcap_{j=1}^{[t/k]-1}
\left\{X_{jk+1, (j+1)k}^{(k)}\neq a_{1,k}\right\}
\, \bigg| \, X_{1, k}^{(k)}=a_{1,k} \right)\, .
\]
The lemma now follows by recursion.
\end{proof}
\begin{lemma}\label{espdr} There exists a positive constant $C$, such
  that for any positive integer $r$ and any positive integer
  $k$, the following inequality holds
\[
\esp\left(\big(D_1^{(k)}\big)^r\right) \le
 r! k^r \left(\frac{1}{\infd}\right)^{kr}\, .
\]
\end{lemma}

\begin{proof}
The result follows immediately from Lemma \ref{upprod1}.
\end{proof}



\section{Proof of Theorem \ref{bootmarkov}}\label{proofbootm}
We can now start the proof of Theorem \ref{bootmarkov}. We first
observe that
\begin{equation}\label{slight}
\frac{\sqrt{R^{(k)*}_{m_k}-1}}{\sigma^{{(k)}_*}}\left(\mu^{(k)*}-
\hat{\mu}^{(k)} \right)
= \frac{\sum_{i=1}^{m_k}{\ZK}^{(k)*}_i}{\sqrt{\var^{*}
\left(\sum_{l=1}^{m_{k}}{\ZK}_{l}^{(k)*}\right)}}\;.
\end{equation}

We want to prove that the right hand side of
\ref{slight} converges in distribution to a standard normal
distribution, when $k \to +\infty$.  By the Lindeberg-Feller Central
Limit Theorem for double arrays (see, for instance,
Billingsley 1999\nocite{bil}),
this will follow once we show that for any $\epsilon >0$
\begin{equation}\label{LFcond}
\lim_{k\to+\infty}
\frac{\esp^*\left( \left({{\ZK}_{1}^{(k)*}}\right)^2\; \one_{\left\{
\left({\ZK}_{1}^{(k)*}\right)^2> \epsilon m_k\var^*\big({\ZK}_{1}^{(k)*}\big)
\right\}}\right)}{\var^* ({\ZK}_{1}^{(k)*} )}= 0\, .
\end{equation}
Using Lemma \ref{expresboot} we can rewrite (\ref{LFcond}) as
\begin{equation}\label{LFbis}
\lim_{k\to+\infty}
\frac{\sum_{l=1}^{m_{k}} {\big({\ZK}_{l}^{(k)}\big)^2\; 
\one_{ \left\{\left({\ZK}_{l}^{(k)}\right)^2 > 
\epsilon \sum_{j=1}^{m_{k}}\big({\ZK}_{j}^{(k)}\big)^{2}\right\}}}}
{\sum_{l=1}^{m_{k}} {\left({\ZK}_{l}^{(k)}\right)}^2}=0\, .
\end{equation}
Since
\begin{equation}
\one_{ \left\{\left({\ZK}_{l}^{(k)}\right)^2 > 
\epsilon \sum_{j=1}^{m_k}\big({\ZK}_{j}^{(k)}\big)^{2}\right\}}
\le
\frac{\left({\ZK}_{l}^{(k)}\right)^2} 
{\epsilon\sum_{l=1}^{m_k}{ \big({\ZK}_{l}^{(k)}\big)^2}}
\, \, ,
\end{equation}
the fraction at the left-hand side of expression \ref{LFbis} is
bounded above by
\begin{equation}\label{frac}
\frac{\sum_{l=1}^{m_k} {\big({\ZK}_{l}^{(k)}\big)^4}}
{\epsilon\left(\sum_{l=1}^{m_k}{
\big({\ZK}_{l}^{(k)}\big)^2}\right)^2}\, \, .
\end{equation}
To prove that expression (\ref{frac}) vanishes as $k$ diverges, we will
obtain a sequence of almost sure upper bounds for its numerator and a
sequence of almost sure lower bounds for its denominator.

\begin{lemma}\label{upperbound} 
For any $\alpha > 0$ and for any $v >1+4\ln(1/\infd)/{\alpha}$, if
 $m_k=\left[e^{\alpha k}\right]$, then for almost all samples the
 upper-bound
\[
\sum_{i=1}^{m_k}\left({\ZK}_i^{(k)}\right)^4 \le {m_k}^{v}\, ,
\]
holds, for all $k$ large enough.
\end{lemma}
\begin{proof} Markov's inequality and Lemmas \ref{basic} and \ref{espdr} 
imply that 
\begin{equation}\label{Markov's}
\pro\left(\sum_{i=1}^{m_k}\left({\ZK}_i^{(k)}\right)^4 > {m_k}^{v}
\right) \le
\frac{\esp\left(\left({\ZK}_1^{(k)}\right)^4\right)}{m_k^{v -1}}\, ,
\end{equation}
\begin{equation}\label{Markov's}
\pro\left(\sum_{i=1}^{m_k}\left({\ZK}_i^{(k)}\right)^4 > {m_k}^{v}
\right) \le
\frac{C k^4}{m_k^{v -1} \infd^{4k}}\, ,
\end{equation}

where $C>0$ does not depend on $k$. Since by hypothesis $\alpha (v-1)
> 4 \ln (1/\infd)$, we conclude that the right hand side of expression
(\ref{Markov's}) is summable. This together with the Borel-Cantelli
Lemma concludes the proof of the lemma.
\end{proof}


The next step is to find a lower bound for the denominator.
\begin{lemma}\label{lowerboundtilde}
For any $\alpha > 4 \ln (1/\infd)$, and for any summable sequence of
non negative real numbers $\eta_k$, $k=1,2,\dots$, if
$m_k=\big[e^{\alpha k}\big]$, then, for almost all samples, the lower
bound
\[
\sum_{i=1}^{m_k}\left({\ZM}_i^{(k)}\right)^2 \ge {m_k}{\eta_k}
\esp\left(\left({\ZM}_1^{(k)}\right)^2\right)  \, ,
\]
holds , for all $k$ large enough.
\end{lemma}
\begin{proof}
To simplify the notation, let us call
\[
W^{(k)}=\sum_{i=1}^{m_k}\left({\ZM}_i^{(k)}\right)^2\, .
\]
By definition we have
\begin{equation}\label{espw}
\esp\left(W^{(k)} \right)=m_k\esp\left(\left({\ZM}_1^{(k)}\right)^2\right).
\end{equation}
Using the fact that the random variables 
\[
\left({\ZM}_i^{(k)}\right)^2-\esp\left(\left({\ZM}_1^{(k)}\right)^2\right)
\]
are independent, identically distributed and have zero mean we get
\begin{equation}\label{espw2}
\esp\left(\left(W^{(k)}\right)^2 \right)=
m_k(m_k-1) \left(\esp\left(\left({\ZM}_1^{(k)}\right)^2\right)\right)^2+
m_k\esp\left(\left({\ZM}_1^{(k)}\right)^4\right)
\, .
\end{equation}
Using the inequality of Paley-Zygmund, for $0<\eta<1$, together with
the identities (\ref{espw}) and (\ref{espw2}) we obtain the inequality
\[
\pro\left(W^{(k)}\ge\eta \esp(W^{(k)})\right) \ge 
\frac{(1-\eta)^{2}m_k^2\left(\esp\left(\left({\ZM}_1^{(k)}\right)^2\right)\right)^2}
{m_k(m_k-1) \left(\esp\left(\left({\ZM}_1^{(k)}\right)^2\right)\right)^2+
m_k\esp\left(\left({\ZM}_1^{(k)}\right)^4\right)}\,.
\]
The right hand-side of the above expression can be rewritten as
\begin{equation}\label{PZ}
(1-\eta)^{2}\left(1-\frac{1}{m_k}+
\frac{\esp\left(\left({\ZM}_1^{(k)}\right)^4\right)}
{{m_k}\left(\esp\left(\left({\ZM}_1^{(k)}\right)^2\right)\right)^2}
\right)^{-1}\, .
\end{equation}
Therefore Lemma \ref{basic} and Hypothesis \ref{Ak} imply that 
\begin{equation}
\pro\left(W^{(k)}\ge\eta \esp(W^{(k)})\right) \ge
(1-\eta)^{2}\left(1-\frac{1}{m_k}+
\frac{C\esp\left(\left(D_1^{(k)}\right)^4\right)}
{m_k}
\right)^{-1}\, ,
\end{equation}
where $C>0$ does not depend on $k$.
From this it follows immediately that
\begin{equation}\label{PZ2}
\pro\left(W^{(k)}<\eta \esp(W^{(k)})\right) \le
\frac{-\frac{1}{m_k}+
\frac{C\esp\left(\left(D_1^{(k)}\right)^4\right)}
{m_k}+2\eta -{\eta}^2}{1-\frac{1}{m_k}+
\frac{C\esp\left(\left(D_1^{(k)}\right)^4\right)}
{m_k}}\, .
\end{equation}
Lemma \ref{espdr} and the choice of $\alpha$ imply that the quantity
\begin{equation}
\left|\frac{1}{m_k}-
\frac{C\esp\left(\left(D_1^{(k)}\right)^4\right)}{m_k}\right|
\le \frac{1}{2}
\end{equation}
for $k$ large enough. Therefore inequality (\ref{PZ2}) implies that
\begin{equation}\label{PZ3}
\pro\left(W^{(k)}<\eta \esp(W^{(k)})\right) \le
2\frac{C\esp\left(\left(D_1^{(k)}\right)^4\right)}
{m_k}+4\eta\, ,
\end{equation}
for $k$ large enough. Using again Lemma \ref{espdr} it follows from
(\ref{PZ3}) that
\begin{equation}
\sum_{k=1}^{+\infty}\pro\left(W^{(k)}<\eta_k \esp(W^{(k)})\right)
<+\infty\, ,
\end{equation}
for any summable sequence of non negative real numbers $\eta_k$,
$k=1,2,\dots$. As a consequence, the Lemma of Borel-Cantelli implies
that
\begin{equation}\label{radenomi}
\sum_{i=1}^{m_k}\left({\ZM}_i^{(k)}\right)^2 \ge \eta_km_k
\esp\left(\left({\ZM}_1^{(k)}\right)^2\right)\, ,
\end{equation}
almost surely for $k$ large enough.
\end{proof}
\begin{lemma}\label{limtilde}
For any $\alpha > 4 \ln (1/\infd)$, if $m_k=\big[e^{\alpha k}\big]$,
then, for almost all samples, the following limit holds
\[
\lim_{k\to+\infty}
\frac{\sum_{l=1}^{m_k} {\big({\ZK}_{l}^{(k)}\big)^4}}
{\left(\sum_{l=1}^{m_k}{
\big({\ZM}_{l}^{(k)}\big)^2}\right)^2} =0\, \, .
\]
\end{lemma}
\begin{proof}
The result follows at once from Lemmas \ref{upperbound} and
\ref{lowerboundtilde} and the Borel-Cantelli Lemma by taking
$1+4\ln(1/\infd)/\alpha < v < 2$ and, for instance, $\eta_k=1/k^2$.
\end{proof}

The expression in the statement of the above lemma is similar to
(\ref{frac}) with ${\ZK}_{l}^{(k)}$ replaced by ${\ZM}_{l}^{(k)}$ in
the denominator. Therefore to conclude
 the proof of Theorem
\ref{bootmarkov} we need the following lemma.

\begin{lemma}\label{dur}
For any $\alpha > 5 \ln (1/\infd)$, if $m_k=\big[e^{\alpha k}\big]$,
then, for almost all samples, the following limit holds
\[
\lim_{k\to+\infty}
\frac{\sum_{l=1}^{m_k}\big({\ZK}_{l}^{(k)}\big)^2}
{\sum_{l=1}^{m_k}\big({\ZM}_{l}^{(k)}\big)^2} =1\, \, .
\]
\end{lemma}

\begin{proof}
An elementary computation shows that for any real numbers $a$ and $b$,
and for any $\epsilon>0$ one has
$$
(1-\epsilon)a^{2}+(1-\epsilon^{-1})b^{2}\le
(a+b)^{2}\le (1+\epsilon)a^{2}+(1+\epsilon^{-1})b^{2}\;.
$$ We apply this inequality for each $l=1,\ldots,m_{k}$ with
$a={\ZM}_{l}^{(k)}$, and $b=(\hat{\mu}^{(k)}-\mu^{(k)}){D}_{l}^{(k)}$.
Summing up over $l$ and using identity (\ref{dif}) we obtain the
inequalities
\begin{equation}\label{verify}
1-\epsilon+(1-\epsilon^{-1})\zeta^{(k)}
\le  \frac{\sum_{l=1}^{m_k}\big({\ZK}_{l}^{(k)}\big)^2}
{\sum_{l=1}^{m_k}\big({\ZM}_{l}^{(k)}\big)^2}\le
1+\epsilon+(1+\epsilon^{-1})\zeta^{(k)}
\end{equation}
where
\begin{equation}\label{zeta}
\zeta^{(k)}=\big(\hat{\mu}^{(k)}-\mu^{(k)}\big)^{2}
\frac{\sum_{l=1}^{m_k}\big({D}_{l}^{(k)}\big)^2}
{\sum_{l=1}^{m_k}\big({\ZM}_{l}^{(k)}\big)^2}\;.
\end{equation}
To conclude the proof it remains to show that $\zeta^{(k)}$ converges to
zero almost surely as $k$ diverges. 

Using Lemma \ref{suchthat}, Markov's inequality and the Borel-Cantelli
Lemma, it follows immediately that for any summable 
 sequence of positive numbers $\rho_{k}$, $k\ge1$, and for almost
all samples, the following inequality holds 

\begin{equation}\label{mumu}
 \big(\hat{\mu}^{(k)}-\mu^{(k)}\big)^{2}\le \frac{C}{\rho_{k}}\left[
\frac{\esp\left(\left(D_1^{(k)}\right)^3\right)}{m_{k}}
+\frac{\esp\left(\left(D_1^{(k)}
\right)^4\right)}{m_{k}^{2}}\right]\;,
\end{equation}
for all $k$ large enough, where $C$ is a positive constant independent
of $k$.
We also observe that for the same sequence $\rho_{k}$ the inequality
\begin{equation}\label{dede}
\sum_{l=1}^{m_k}\big({D}_{l}^{(k)}\big)^2\le \frac{m_{k}}{\rho_{k}}
\esp\left(\left(D_1^{(k)}
\right)^2\right)\;,
\end{equation}
holds almost surely for all $k$ large enough.

Combining Lemma \ref{lowerboundtilde} and Hypothesis \ref{Ak}, we conclude
that for any summable sequence $\eta_{k}$, $k\ge 1$, and for almost all
sample, the following inequality holds  
\begin{equation}\label{enfin}
\sum_{l=1}^{m_k}\big({\ZM}_{l}^{(k)}\big)^2\ge Cm_{k}\eta_{k}\;,
\end{equation}
for all $k$ large enough, 
where $C$ is a strictly positive constant independent of $k$.

Using  inequalities (\ref{mumu}), (\ref{dede}), (\ref{enfin}), and using
Lemma \ref{espdr} we deduce  that for almost all samples, the following
inequality holds 
$$
\zeta^{(k)}\le 
C\frac{e^{-k(\alpha-5\ln(1/\infd))}}{\rho_{k}^{2}\eta_{k}}
$$
for all $k$ large enough, where $C$ is a positive constant independent of
$k$. Since by hypothesis, $\alpha>5\ln(1/\infd)$, it is enough to take
for instance $\rho_{k}=\eta_{k}=1/k^{2}$ to conclude $\zeta^{(k)}$
converges to zero almost surely. Recalling that inequality
(\ref{verify}) holds for any fixed $\epsilon>0$, the lemma follows.
\end{proof}

Combining Lemmas \ref{limtilde} and \ref{dur}, it follows that almost
surely 
\begin{equation}\label{encore}
\lim_{k\to+\infty}\frac{\sum_{l=1}^{m_k} {\big({\ZK}_{l}^{(k)}\big)^4}}
{\epsilon\left(\sum_{l=1}^{m_k}{
\big({\ZK}_{l}^{(k)}\big)^2}\right)^2}=0\, \, .
\end{equation}
This implies  (\ref{LFcond}) and finishes the proof of Theorem
\ref{bootmarkov}. 
\section{Proof of Theorem \ref{boot}.}\label{proofboot}

The basic idea of the proof is to approximate the chain of infinite
order by a sequence of Markov chains of increasing order satisfying
the hypotheses of Theorem \ref{bootmarkov}.  We will use for this
purpose the canonical Markov approximation $\XC$ of the chain $\X$
which is the Markov chain of order $k$ whose transition probabilities
are defined by
\begin{equation}
\label{defapproxMarkov}
P^{[k]}(b \; | \; a_1, \ldots, a_k) := \pro(X_{k+1} = b | X_j = a_j, 1
 \leq j \leq k)
\end{equation}
for all integer $k\geq 1$ and $a_1, \ldots, a_k, b \in A$.  

From now on we only consider stationary chains.  The sequence of
stationary canonical Markov approximations can be constructed together
with the stationary chain of infinite order on the same probability
space $(\Omega, \mathcal{F}, \pro)$. In particular they can be
constructed together using the well-known {\sl maximal coupling}(see,
for instance, Appendix A.1 in Barbour Holst and Janson,
1992\nocite{barholjan92}). For details of this construction in the
present context we refer the reader to Fern\'andez and Galves
(2002)\nocite{fg}.

Before starting the proof of Theorem \ref{boot} we will recall a few
results from the literature which will be used in the sequel. The
following theorem was proven in Fern\'andez and Galves
(2002)\nocite{fg}.

\medskip
\noindent{\bf Theorem.} {\sl Let $\X$ be a chain of infinite order on
the finite alphabet $A$ and
satisfying the conditions
\[
\sum_{a \in A}\inf_{ (\ldots, x_{-2}, x_{-1}) \in  {\cal{A}}_a } p(a
|x_{-1}, x_{-2}, \ldots) >0
\quad\mathrm{and}\quad
\sum_{l\ge 1}\beta_l <+\infty\; .
\]
Then the construction of the chains using the maximal coupling
satisfies the following inequality
\begin{equation}\label{probdif}
\pro\left\{X^{[k]}_0\neq X_0 \right\}\le \beta_k\; .
\end{equation}}

The following theorem is a particular case of the main theorem of
Bressaud, Fern\'andez and Galves (1999)\nocite{bfg2}. For convenience
of the reader we will reformulate the result in the framework in which
it will be used in the proofs below.
\medskip

\noindent{\bf Theorem.} {\sl If hypotheses ${\bf H_1}$
and ${\bf H_2}$ are satisfied then the chain $\X$ is exponentially
$\varphi$-mixing}.

\medskip
For a definition of $\varphi$-mixing chains we refer the reader to
Billingsley (1999)\nocite{bil}. To make the connection between the
present hypotheses and the assumptions of Bressaud {et al.}
(1999)\nocite{bfg2} we note that hypotheses ${\bf H_1}$ and ${\bf
H_2}$ imply that the sequence of log-continuity rates $(\gamma_l)$
defined by
\[
\gamma_l=\max_{a \in A}\; 
\sup_{ \stackrel
{(\ldots, x_{-2}, x_{-1}) \in  {\cal{A}}_a} 
{{x_i=y_i}\; , \; {i=-l,\ldots, -1}}}\; 
\left|\frac{p(a |x_{-1},x_{-2},\ldots)}{p(a |
y_{-1},y_{-2},\ldots)}-1\right|
\]
is exponentially decreasing and therefore satisfies the hypotheses of
this paper.
\medskip

We can now start the proof of Theorem \ref{boot}. First of all we will
use the above mentioned result by Fern\'andez and Galves (1999) to
obtain an upper bound for the probability of discrepancies in the
first $r$ symbols for the coupled realizations of the chain $\X$ and
its canonical Markov approximation of order k $\XC$. More precisely
let us define
\[
\Delta_r^{[k]} \,:=\, \{X_t^{[k]} = X_t, \, t = 1 \ldots,
r \}\; ,
\]
which is the set of coincidence up to time $r$ of the chains $\XC$ and
$\X$.

\begin{lemma}  \label{lfg} Let  $\X$ be  a chain of infinite order 
satisfying conditions ${\bf H_1}$ and ${\bf H_2}$  with ${\beta}_l$
summable. The there exists a positive constant $C$ such that
\[
\pro\left\{\left(\Delta_r^{[k]}\right)^c\right\} \le C r \beta_k
\]
\end{lemma}

We will now check that the hypotheses of Theorem \ref{bootmarkov} are
satisfied by the sequence of canonical Markov approximations $\XC$, $k
\ge 1$.

\begin{lemma}\label{deltabarposi} Under assumption ${\bf H_1}$ we have
\[
 \inf\{\delta^{[k]}:\; k\ge 1\} \ge \delta\; ,
\]
where
\[
\delta^{[k]}= 
\min_{a \in A}\inf_{ ( x_{-k}, \ldots , x_{-1})\in {\cal{A}}^{(k)}_a } 
p^{[k]}(a |x_{-1}, \ldots, x_{-k})\; .
\]
\end{lemma}
\begin{proof} 
Follows at once from the properties of the conditional probability.
\end{proof}
This lemma establishes condition (\ref{deltauni}). The proof that
condition (\ref{Ak}) holds follows from the next three lemmas. Let us
define
\[
Z_i(k)=
\sum_{n=R_{i-1}(k-1)}^{R_{i}(k)-1}
\left(f(X_n)-\mu\right)
\quad \mathrm{and} \quad
Z^{[k]}_i=
\sum_{n=R^{[k]}_{i-1}}^{R^{[k]}_{i}-1}
\left(f\left(X^{[k]}_n\right)-\mu^{[k]}\right)
\; ,
\]
where $R_1^{[k]}$ is defined as in expression (\ref{recursi}) 
using the chain $\XC$ and $\mu^{[k]}= \esp\big(
f(X^{[k]}_1)\big)$. 

\medskip
\begin{lemma}\label{morceau} Under Hypotheses ${\bf
H_1}$, ${\bf H_2}$ and ${\bf H_3}$ the chain $\X$ satisfies the
inequality
\[
\liminf_{k\to+\infty}\esp\left(\left(Z_1(k)\right)^2\right) >0\; .
\]
\end{lemma}
\begin{proof}
Markov's inequality implies that
\[
\esp\left(\left(Z_1(k)\right)^2\right)
 \ge 
u^2 \pro\left\{\left(Z_1(k)\right)^2 >u^2 \right\}\; ,
\]
for any real number $u$.
Recalling that $R_{1}(k)\ge 1$, we obtain the lower bound
\begin{equation}\label{borneinf}
\esp\left(\left(Z_1(k)\right)^2\right)
\ge
u^2\pro\left\{\frac{\left|Z_1(k)\right|}{\sqrt{R_{1}(k)}} >u \right\}
\; .
\end{equation}
By the above mentioned theorem from Bressaud {et
al.} (1999)\nocite{bfg2}, the process $(f(X_n))_n$ is exponentially
$\varphi$-mixing. Therefore it follows from classical results on the
Central Limit Theorem (cf. for instance Theorems 20.1 and 20.3 from
Billingsley 1999\nocite{bil})
\[
\frac
{ Z_1(k)}
{\sqrt{R_{1}(k)}}\stackrel{{\mathcal D}}
{\longrightarrow} {\cal N}(0,\sigma^2)
\]
as $k$ diverges. Hypothesis ${\bf H_3}$ ensures that $\sigma>0$.  This
implies that for any fixed $u$ and any $k$ large enough the lower bound
provided by inequality (\ref{borneinf}) is greater than a fixed
strictly positive real number. This concludes the proof of the lemma.
\end{proof}

We define $D_{l}(k)=R_{l}(k) -R_{l-1}(k)$.

\begin{lemma}\label{upprod2}
For any integer $k$, any integer $r \le 4$ and any positive real
number $t$ the following inequalities hold
\[
\pro\left(D_1(k) >t\right) \le
\left(1-\delta^k\right)^{[t/k]}
\; \quad \mathrm{and}\quad\;
\esp\left(\big(D_1(k) \big)^r\right) \le
 C k^r \left(\frac{1}{\delta}\right)^{kr}\, .
\]
where $C$ is a positive constant.
\end{lemma}

\begin{proof} 
The proof is exactly the same as the proofs of Lemmas \ref{upprod1}
and \ref{espdr}. 
\end{proof}

\begin{lemma}\label{saispas}
Under the conditions of Theorem \ref{boot} the sequence of canonical
Markov approximations satisfies the inequality
\[
\liminf_{k\to+\infty}\esp\left(\left(Z^{[k]}_1\right)^2\right)>0\; .
\]
\end{lemma}

\begin{proof}
We will first derive an upper bound for the the modulus of the difference
\[
\left|\esp\left(\big(Z_1(k)\big)^2-
\big(Z^{[k]}_1\big)^2\right)\right|=
\left|\esp\left((Z_1(k)-Z^{[k]}_1)(Z_1(k)+Z^{[k]}_1) \right)\right|
\]
The finiteness of the alphabet $A$ implies that
\begin{equation}\label{finitude}
\left|Z_1(k)+Z^{[k]}_1\right|\le C\left|R_1(k)+R^{[k]}_1 \right|\, ,
\end{equation}
where $C=\max\{|f(a)|: a\in A\}$. 
We observe also that
\begin{equation}\label{comb2}
\left|Z_1(k)-Z^{[k]}_1\right|
\le
\sum_{n=1}^{R_1(k)\wedge R^{[k]}_1}\left|Y_n-Y^{[k]}_n\right|
+C\left|R_1(k)- R^{[k]}_1\right|\; ,
\end{equation}
where $ Y_n=f\left(X_n\right)-\mu$ and
$Y^{[k]}_n=f\left(X^{[k]}_n\right) -\mu^{[k]}$.

In the sequel we will no longer specify the different positive
constants appearing in the various estimates. Moreover they will be
all denoted by the letter $C$.  Combining inequalities
(\ref{finitude}) and (\ref{comb2}) we obtain
\begin{equation}\label{endeux}
\left|\esp\left(\big(Z_1(k)\big)^2-
\big(Z^{[k]}_1\big)^2\right)\right|\le \Oun \esp\left(
\big|R_1(k)- R^{[k]}_1\big|\big|R_1(k)+ R^{[k]}_1\big|\right)
$$
$$
+\Oun \esp\left(\sum_{n=1}^{R_1(k)\wedge
R^{[k]}_1}\left|Y_n-Y^{[k]}_n\right|
\big|R_1(k)+ R^{[k]}_1\big|\right)\;.
\end{equation}

We will estimate separately each term.
For the second term we have
$$
\esp\left(\sum_{n=1}^{R_1(k)\wedge
R^{[k]}_1}\left|Y_n-Y^{[k]}_n\right|
\big|R_1(k)+ R^{[k]}_1\big|\right)
$$
$$
= \esp\left(\gun_{(\Delta_{k}^{[k]})^{c}}\sum_{n=1}^{R_1(k)\wedge
R^{[k]}_1}\left|Y_n-Y^{[k]}_n\right|
\big|R_1(k)+ R^{[k]}_1\big|\right)
$$
$$
\le \Oun
 \esp\left(\gun_{(\Delta_{k}^{[k]})^{c}}
\big(R_1(k)+ R^{[k]}_1\big)^{2}\right)\;.
$$
Using Schwarz inequality and Lemmas \ref{espdr}, \ref{lfg} and
\ref{upprod2}. we obtain the upper bound
$$
\esp\left(\gun_{(\Delta_{k}^{[k]})^{c}}\right)^{1/2}
 \esp\left(\big(R_1(k)+ R^{[k]}_1\big)^{4}\right)^{1/2}
\le
\Oun k^{5/2}\beta_{k}^{1/2}\delta^{-2k}\; .
$$

We now come to the estimation of the first term in (\ref{endeux}).
Using Scwharz inequality and Lemmas \ref{espdr} and \ref{upprod2} we
get
$$
 \esp\left(
\big|R_1(k)- R^{[k]}_1\big|\big|R_1(k)+ R^{[k]}_1\big|\right)
\le  \esp\left(
\big(R_1(k)- R^{[k]}_1\big)^{2}\right)^{1/2}
 \esp\left(\big(R_1(k)+ R^{[k]}_1\big)\right)^{1/2}
$$
$$
\le\Oun k \delta^{-k}\esp\left(
\big(R_1(k)- R^{[k]}_1\big)^{2}\right)^{1/2}\;.
$$
We now have
$$
\esp\left(
\big(R_1(k)- R^{[k]}_1\big)^{2}\right)=
\esp\left(\gun_{\Delta_{k}^{[k]}}
\big(R_1(k)- R^{[k]}_1\big)^{2}\right)+\esp\left(\gun_{(\Delta_{k}^{[k]})^{c}}
\big(R_1(k)- R^{[k]}_1\big)^{2}\right)
$$
and the last term is estimated as above. For the first term, we have
$$
\esp\left(\gun_{\Delta_{k}^{[k]}}
\big(R_1(k)- R^{[k]}_1\big)^{2}\right)\le
\esp\left(\big(R_1(k)+ R^{[k]}_1\big)^{2}\left(1-\prod_{j=R_1(k)\wedge
R^{[k]}_1}^{R_1(k)\wedge R^{[k]}_1+k-1}\gun_{X_{j}^{[k]}=X_{j}}\right)\right)
$$
$$
\le \esp\left(\big(R_1(k)+ R^{[k]}_1\big)^{4}\right)^{1/2}
 \esp\left(\left(1-\prod_{j=R_1(k)\wedge
R^{[k]}_1}^{R_1(k)\wedge
R^{[k]}_1+k-1}\gun_{X_{j}^{[k]}=X_{j}}\right)\right)^{1/2} 
$$
$$
\le \Oun k^{2}  \delta^{-2k}\esp\left(\left(1-\prod_{j=R_1(k)\wedge
R^{[k]}_1}^{R_1(k)\wedge
R^{[k]}_1+k-1}\gun_{X_{j}^{[k]}=X_{j}}\right)\right)^{1/2} 
$$
where we have used again Schwarz inequality and Lemmas \ref{espdr} and
\ref{upprod2}.
We now have
\[
\esp\left(\left(1-\prod_{j=R_1(k)\wedge
R^{[k]}_1}^{R_1(k)\wedge
R^{[k]}_1+k-1}\gun_{X_{j}^{[k]}=X_{j}}\right)\right)
=\sum_{p=1}^{\infty}\esp\left(\gun_{R_1(k)\wedge
R^{[k]}_1=p}\left(1-\prod_{j=p}^{p+k-1}
\gun_{X_{j}^{[k]}=X_{j}}\right)\right)\; .
\]
Using Schwarz inequality and stationarity and Lemmas \ref{espdr},
\ref{lfg} and \ref{upprod2}.  this is bounded above by
$$
\left(\sum_{p=1}^{\infty}p^{2} \esp\left(\gun_{R_1(k)\wedge
R^{[k]}_1=p}\right)\right)^{1/2}
\esp\left(\gun_{(\Delta_{k}^{[k]})^{c}}\right)^{1/2} 
$$
$$
\le \esp\left(\big(R_1(k)+R^{[k]}_1\big)^{2}\right)^{1/2}
\esp\left(\gun_{(\Delta_{k}^{[k]})^{c}}\right)^{1/2} 
\le\Oun k^{3/2} \delta^{-k}\beta_{k}^{1/2} 
$$
Collecting together the above bounds we get
$$
\left|\esp\left(\big(Z_1(k)\big)^2-
\big(Z^{[k]}_1\big)^2\right)\right|
\le \Oun\left(k^{5/2} \delta^{-2k}\beta_{k}^{1/2} +k^{19/8}
 \delta^{-9k/4}\beta_{k}^{1/8}\right)\;.
$$ 
It follows from this inequality and assumption $c>18\log \delta^{-1}$ that
$$
\lim_{k\to\infty}\left|\esp\left(\big(Z_1(k)\big)^2-
\big(Z^{[k]}_1\big)^2\right)\right|=0\;.
$$
This together with Lemma \ref{morceau} concludes the proof of the lemma.
\end{proof}

In order to prove Theorem \ref{boot} we need to construct together the
bootstrap samples of $\X$ and $\XC$. We recall that we have already
assumed that $\X$ and $\XC$ are constructed together using the maximal
coupling. Now, given two coupled realizations of theses chains we will
use the same realization of the sequence of random indices to choose
the blocks entering in the bootstrap samples of the chains.  Formally,
for every fixed $k \ge 1$ the bootstrap blocks will be defined as
\[
\xi^{\ast}_{l}(k)=\xi_{I_l(k)}(k)\, \quad \, \mathrm{and}\quad\, 
\xi^{[k]\ast}_{l}=\xi^{[k]}_{I_l(k)}
\]
where $I_1(k),\ldots,I_{m_k}(k)$ are the same independent random
variables with uniform distribution in the set $\{1,\ldots,m_k\}$.

The next lemma says that the coupled samples of $\X$ and $\XC$
coincide up to time $R_{m_k}(k)$ with overwhelming probability.

\begin{lemma} \label{ffin} 
Under the hypotheses of Theorem \ref{boot} we have
\[
\lim_{k\to +\infty}
\pro\left(\left(\Delta_{R_{m_k}(k)}\right)^c\right)=0\; .
\]
\end{lemma}

\begin{proof}
We observe that for any $r >0$ we have
\begin{equation}\label{ouf1}
\pro\left(\left(\Delta_{R_{m_k}(k)}\right)^c\right)
\le 
\pro\left(\left(\Delta_{r}\right)^c\right)
+
\pro\left(R_{m_k}(k)>r\right)\; .
\end{equation}
By Lemma \ref{lfg} the first term in the right hand side of
(\ref{ouf1}) is bounded above by $C r\beta_k$. 

It follows from Lemmas \ref{upprod2} and \ref{deltabarposi} that the
second term of the right hand side of (\ref{ouf1}) is bounded above by
\[
m_k\pro\left(D_{1}(k)>r/m_k\right)
\le
m_k \left(1-\delta^k\right)^{[r/(k m_k)]}
\;.
\]
We now set $r=\lambda k^2 m_k \delta^{-k}$, where $\lambda$ is a fixed
number strictly larger than $\alpha$.  With this choice of $r$ the two
terms in inequality (\ref{ouf1}) tends to $0$ when $k$ diverges. This
concludes the proof of the lemma.
\end{proof}

We can now conclude the proof of Theorem \ref{boot}. First of all we
observe that
\[
\frac{\sqrt{R^*_{m_k}(k)-1}}{\sigma^*_k}
 \left(\mu^{\ast}_k - \hat{\mu}_k\right)= 
\frac{\sigma^{[k]*}}{\sigma^*_k}
\sqrt{\frac{R^*_{m_k}(k)-1}{R^{[k]*}_{m_k}-1}}\quad
\frac{\sqrt{R^{[k]\ast}_{m_k}-1}}{\sigma^{[k]\ast}}
\left(\mu^{[k]\ast} - \hat{\mu}^{[k]}  \right)
\]
\[
+\frac{\sqrt{R^*_{m_k}(k)-1}}{\sigma^*_k}
\left (\hat{\mu}^{[k]}-\hat{\mu}_{k} \right)
+
\frac{\sqrt{R^*_{m_k}(k)-1}}{\sigma^*_k}\left ( \mu^*_{k} -
\mu^{[k]*}\right )
\]

Lemma (\ref{ffin}) ensures that last two terms are equal to zero 
with probability tending to 1 when
$k$ tends to infinity. Theorem \ref{bootmarkov} implies 
\[
\sqrt{\frac{R^*_{m_k}(k)-1}{R^{[k]*}_{m_k}-1}}\quad
\frac{\sqrt{R^{[k]\ast}_{m_k}-1}}{\sigma^{[k]\ast}}
\left(\mu^{[k]\ast} - \hat{\mu}^{[k]}  \right)
\stackrel{{\mathcal D}}{\longrightarrow}
{\cal N}(0,1)\; .
\]
Finally we observe that Lemma (\ref{ffin})  ensures that 
\[
\lim_{k\to\infty}\pro\left(\frac{\sigma^{[k]*}}{\sigma^*_k}
\sqrt{\frac{R^*_{m_k}(k)-1}{R^{[k]*}_{m_k}-1}}=1\right)=1\;.
\]
This concludes the proof of Theorem \ref{boot}. 

\section*{Acknowledgments} We thank M.~Cassandro, R.~Fern\'andez, 
D.~Gabrielli and N.~Garcia for helpful discussions.


\end{document}